\documentclass[11pt]{amsart}
\usepackage[T1]{fontenc}
\usepackage[utf8]{inputenc}
\usepackage{lmodern,amsmath,amssymb,amsthm}
\usepackage[margin=1in]{geometry}
\usepackage[hidelinks]{hyperref}
\newtheorem{theorem}{Theorem}[section]
\newtheorem{proposition}[theorem]{Proposition}
\newtheorem{lemma}[theorem]{Lemma}

\theoremstyle{definition}
\newtheorem{definition}[theorem]{Definition}
\newtheorem{remark}[theorem]{Remark}
\newcommand{\Q}{\mathbb Q}
\newcommand{\Z}{\mathbb Z}
\newcommand{\F}{\mathbb F}
\newcommand{\Jac}{\operatorname{Jac}}
\newcommand{\Disc}{\operatorname{Disc}}

\newcommand{\ceil}[1]{\left\lceil #1\right\rceil}
\newcommand{\floor}[1]{\left\lfloor #1\right\rfloor}
\title[Arithmetic circuits for factorial multiples]{Unconditional Subexponential Arithmetic Circuits\ for Positive Factorial Multiples}
\author{Liang Zhang}
\email{tensor@pku.edu.cn}
\thanks{Corresponding author: Liang Zhang.}
\date{}

\begin{document}
\begin{abstract}
Let $\tau(n)$ be the least number of binary additions, subtractions and
multiplications required to construct $n$ from $1$ and $2$, with reuse of
intermediate values. We show that, for every sufficiently large integer $k$,
some positive integer $m$ satisfies
\[
 \tau(mk!)\leq
 \exp\!\left(O\!\left((\log k)^{2/3}(\log\log k)^{1/3}\right)\right).
\]
The proof combines unconditional results of Lenstra, Pila and Pomerance
on smooth numbers in short intervals and on orders of genus-two Jacobians
with a five-parameter family carrying a rational divisor class of infinite order.
We encode rational divisor arithmetic, including the tests determining its
branches, by a nonzero division-free integer circuit. A common finite
sampling space yields circuits whose outputs collectively contain every
prime at most $k$, and repeated squaring supplies the required factorial
valuations. The circuit bound includes the cost of constructing all integer
parameters. This work was carried out with assistance from the UniMath model,
developed by UniDT (Shanghai) Co., Ltd., and GPT-6.1 Sol.
\end{abstract}
\maketitle
\setcounter{tocdepth}{1}
\tableofcontents

\section{Introduction}
How many arithmetic operations are needed to construct a given integer?
This question provides a concrete setting in which to study the power and
limitations of arithmetic computation. Although repeated squaring produces
integers of enormous size using very few operations, counting arguments
show that most integers admit no comparably short description. For explicit
sequences, such as the factorials, the gap between these two phenomena
leads to difficult questions connected with integer factorization,
arithmetic circuit lower bounds and algebraic versions of the $P$ versus
$NP$ problem.

Following de Melo and Svaiter~\cite{dMS}, an integer straight-line program
starts with $1,2$ and uses binary $+,-,\times$ gates, with unrestricted
reuse of intermediate values. Its length counts gates, irrespective of
operand bit length. We use this convention throughout. Starting with $1$
alone changes the minimum length by at most one, so it does not affect
the asymptotic questions considered here. Divisions, arbitrary constants,
modular operations, determinants and group operations are not additional
gates.

\begin{definition}
For an integer $n$, let $\tau(n)$ be its minimum straight-line program length.
For $k\geq2$ put
\[
 \mu(k)=\min_{m\in\Z_{>0}}\tau(mk!).
\]
\end{definition}

The function $\mu$ measures the \emph{ultimate complexity} of factorials:
one may choose a convenient multiple of $k!$, but that multiple must still
be constructed from the prescribed constants. Allowing negative nonzero
multipliers gives an equivalent complexity measure up to an absolute
additive constant, since a sign can be changed with at most two gates.

In \emph{The cost of computing integers}, de Melo and
Svaiter~\cite{dMS} established the elementary bounds
\[
 \log_2\log_2 n\leq\tau(n)\leq2\log_2 n\qquad(n\geq2)
\]
and proved that, for every $\varepsilon>0$, almost all integers require
at least order $\log n/(\log\log n)^{1+\varepsilon}$ operations.
Here ``almost all'' means that the exceptional set has natural density
zero. This answered negatively Smale's question whether every integer
could be constructed in a number of operations polynomial in its
double logarithm. Moreira~\cite{Moreira} subsequently obtained sharper
estimates: up to the immaterial choice of initial constants,
\[
 \begin{aligned}
 \tau(n)&\geq(1-o(1))\frac{\log n}{\log\log n}
 &&\text{for almost all }n,\\
 \tau(n)&\leq(1+o(1))\frac{\log n}{\log\log n}
 &&\text{for all sufficiently large }n.
 \end{aligned}
\]
These results describe the general scale of integer complexity, but
density statements do not give lower bounds for a specified sparse
sequence, or for the easiest multiples of its terms.

The factorial questions recorded in~\cite[p.~1378]{dMS} ask whether
$\tau(k!)$, and more generally $\mu(k)$, are bounded by a polynomial in
$\log k$. Their significance comes from the work of Shub and
Smale~\cite{SS}: if $\mu(k)$ is not polynomially bounded in $\log k$,
then Hilbert's Nullstellensatz is intractable in the Blum--Shub--Smale
model over $\mathbb C$, and hence
\[
 P_{\mathbb C}\ne NP_{\mathbb C}.
\]
This is a statement about computation over the complex numbers, rather
than the classical Boolean $P$ versus $NP$ problem. A related route to
the same separation is the Shub--Smale $\tau$-conjecture, which bounds the
number of distinct integer zeros of a nonzero polynomial
$f\in\Z[X]$ by a fixed power of $1+\tau(f)$; here $\tau(f)$ is the
division-free straight-line complexity of the polynomial starting from
$1$ and $X$. Smale included this conjecture as Problem~4 in his
list of mathematical problems for the next century~\cite{Smale}.

Subsequent work clarified why strong lower bounds for these explicit
integer sequences would be consequential. Koiran~\cite{Koiran} showed
that failure of a polylogarithmic bound for $\tau(k!)$ would imply
$VP^0\ne VNP^0$ in the constant-free Valiant model or $P\ne PSPACE$.
B\"urgisser~\cite{Burgisser} proved that if the permanent has
polynomial-size constant-free, division-free arithmetic circuits,
then $\tau(k!)$ is polynomially bounded in $\log k$.
Thus a superpolynomial lower bound in $\log k$ for the factorials
would yield a lower bound for the permanent in that circuit model.
More recently, Alekseev, Grigoriev, Hirsch and Tzameret~\cite{AGHT}
used the Shub--Smale hypothesis on the hardness of factorial multiples
to obtain conditional lower bounds for the Ideal Proof System.
These are implications from hardness assumptions, not established
superpolynomial lower bounds for $\mu(k)$.

On the upper-bound side, successive multiplication immediately gives
$\mu(k)\leq\tau(k!)=O(k)$. Classical fast factorial algorithms, going
back to Strassen, improve the arithmetic-operation count to
$O(\sqrt{k}(\log k)^2)$; see the discussion in
Cheng~\cite[\S1]{Cheng} and Markstr\"om~\cite[\S2]{Markstrom}.
Cheng~\cite[Conjecture~1 and Theorem~1]{Cheng} obtained a substantially
smaller bound for factorial multiples,
\[
 \mu(k)\leq
 \exp\!\left(O\!\left(\sqrt{\log k\log\log k}\right)\right),
\]
assuming a conjecture on the density of smooth integers in elliptic
Hasse intervals. His proof uses elliptic-curve factorization methods
and also provides a randomized subexponential-time construction of the
straight-line program. Complementing these asymptotic results,
Markstr\"om~\cite{Markstrom} used exhaustive searches to determine exact
complexities and upper and lower bounds for small factorials, primorials
and their multiples.

The purpose of this paper is to obtain an unconditional subexponential
upper bound for $\mu(k)$. We use the short-interval smooth-number theorem
and genus-two Jacobian-order estimates of Lenstra, Pila and
Pomerance~\cite{LPP1,LPP2}. Our main result is the following.

\begin{theorem}\label{thm:main}
There is an absolute constant $C>0$ such that for every integer $k\geq2$
there exists a positive integer $m$ satisfying
\begin{equation}\label{eq:allk}
 \tau(mk!)\leq
 \exp\!\left(C[\log(16k^2)]^{2/3}
                    [\log\log(16k^2)]^{1/3}\right).
\end{equation}
Consequently, as $k\to\infty$,
\begin{equation}\label{eq:main}
 \mu(k)\leq
 \exp\!\left(O\!\left((\log k)^{2/3}(\log\log k)^{1/3}\right)\right)
 =k^{o(1)}.
\end{equation}
All logarithms without a subscript are natural.
\end{theorem}

In particular, for every fixed $\varepsilon>0$, one has
$\mu(k)\leq k^\varepsilon$ for all sufficiently large $k$. The exponent
in our upper bound is larger than in Cheng's conditional estimate;
the conclusion here is unconditional. The result concerns multiples
of $k!$, rather than an improved bound for computing $k!$ itself.
It is compatible with the possibility that $\mu(k)$ exceeds every
fixed power of $\log k$, and therefore does not settle Smale's
factorial-multiple question or the associated complexity separations.
The elementary size bound above gives only
$\mu(k)\geq\log_2\log_2(k!)=\Omega(\log k)$, leaving a substantial
gap between the lower and upper bounds.

We briefly describe the proof and its relation to earlier work.
The use of factorization-style arithmetic to obtain short programs
for factorial multiples is already central to Cheng's construction.
The genus-two setting permits unconditional analytic input: the
Jacobian orders supplied by~\cite{LPP2} lie in an interval centered
at $p^2$ of length $p^{3/2}$. With $x$ of order $p^2$, this length
is of order $x^{3/4}$, within the range $y\geq x^{3/5}$ of the
smooth-number theorem in~\cite{LPP1}. An elliptic Hasse interval
instead has length of order $p^{1/2}$ around $p$, outside that range.

To turn these results into integer circuits, we construct a family
of genus-two curves with a marked rational divisor class of infinite
order. Complete Cantor arithmetic computes its scalar multiples using
a bounded number of rational operations per group operation. We
record the nonzero tests specifying the rational branch path, and
encode modular division failures by prime divisors of a nonzero
integer. A counting argument transfers the Jacobian-order estimates
to the marked family. Finally, a common finite sampling space and
a union bound give finitely many certificates covering every prime
at most $k$; repeated squaring supplies its multiplicity in $k!$.
All parameter constructions and arithmetic operations are included
in the circuit length. The theorem asserts existence of such circuits
for each $k$; no comparable bit-complexity bound for finding the
successful samples or their rational branch paths is claimed.

The paper is organized as follows. Section~\ref{sec:inputs} records
the external results. Section~\ref{sec:circuits} gives two elementary
circuit lemmas, and Section~\ref{sec:family} constructs the marked
family. Section~\ref{sec:arithmetic} produces the integer certificates.
Section~\ref{sec:density} proves the smooth-order probability bound,
Section~\ref{sec:cover} obtains simultaneous prime coverage, and
Section~\ref{sec:asymptotics} completes the complexity estimate.

\section{Published inputs}\label{sec:inputs}
We state the two analytic and counting inputs in the exact forms used.
An integer is $B$-smooth if every prime dividing it is at most $B$.

\begin{theorem}[Lenstra--Pila--Pomerance, 1993]\label{thm:lpp1}
Let $c_*=(1980000)^{1/3}$. There is an effectively computable $x_*$ such
that if $x\geq x_*$ and $x^{3/5}\leq y\leq x$, the interval $[x,x+y]$
contains at least
\begin{equation}\label{eq:smooth-count}
 y\exp\!\left(- (\log x)^{1/3}(\log\log x)^{8/3}\right)
\end{equation}
integers that are
$\exp(c_*(\log x)^{2/3}(\log\log x)^{1/3})$-smooth.
\end{theorem}
This is~\cite[Theorem~1.3, p.~399]{LPP1}; the exponent in the density loss
is $8/3$. It must be distinguished from the $1/3$ in the smoothness cutoff.

\begin{theorem}[Lenstra--Pila--Pomerance, 2002]\label{thm:lpp2}
For every prime $p>8100$, all but at most $p$ integers $n$ in
\begin{equation}\label{eq:interval}
 I_p=[p^2-\tfrac12p^{3/2},\ p^2+\tfrac12p^{3/2}]
\end{equation}
occur as $\#\Jac(C_f)(\F_p)$ for at least
\begin{equation}\label{eq:model-count}
 \frac{p^{9/2}}{24000(\log p)^2(\log\log p)^2}
\end{equation}
squarefree quintics $f\in\F_p[X]$, where $C_f$ is the smooth projective
curve with function field $\F_p(X,Y)$, $Y^2=f(X)$.
\end{theorem}
Here quintics have arbitrary nonzero leading coefficient. This is
\cite[Theorem~1.1, p.~106]{LPP2}, a result at each fixed prime, not an
average over primes. It counts polynomial models, not only isomorphism classes.

We also use standard facts about curves and their Jacobians, recorded here
to specify the algebraic prerequisites. A squarefree odd quintic over a field
of characteristic different from two has a smooth genus-two completion
with one rational point at infinity. Its Jacobian is the group of
degree-zero divisor classes. The Hasse--Weil bound is
\begin{equation}\label{eq:weil}
 \big|\#C_f(\F_p)-(p+1)\big|\leq4\sqrt p.
\end{equation}
See~\cite[\S11]{Milne} and~\cite[p.~105]{LPP2} for the Weil-theoretic
setting. Smooth proper curves over a discrete valuation ring have relative
Jacobians that are abelian schemes, compatible with the fibers and the
Abel--Jacobi map~\cite[\S8, Theorem~8.1 and Corollary~12.3]{Milne}.
For an abelian scheme of dimension two over $\Z_3$, formal coordinates at
the identity identify the kernel of reduction on $\Z_3$-points with
$(3\Z_3)^2$, with an integral convergent formal group law. For the
genus-two Jacobians used here, this identification and the integral formal
group law are given explicitly in~\cite[Corollary~4.5 and
Theorem~4.6]{Grant}. The torsion-freeness needed here is proved in
Section~\ref{sec:family}.

The reduced Mumford representation and complete Cantor arithmetic are used
as algebraic theorems, with their formulas displayed in
Section~\ref{sec:arithmetic}. Precise references are
Galbraith~\cite[Theorems~10.3.14, 10.3.21 and 10.4.1]{Galbraith}
(online version~2.0), and Kohel--Smith~\cite[\S2, Algorithm~1]{KS}.
The original algorithm is due to Cantor~\cite{Cantor}.

\section{Two elementary circuit lemmas}\label{sec:circuits}
\begin{lemma}[Prime support to factorial valuations]\label{lem:amplify}
If $k\geq2$ and $q\in\Z\setminus\{0\}$ is divisible by every prime at
most $k$, then appending $r=\ceil{\log_2 k}$ squaring gates produces a
positive multiple of $k!$. Its length is at most $\tau(q)+r$.
\end{lemma}
\begin{proof}
The output is $q^{2^r}>0$. For $p\leq k$,
\[
 v_p(k!)=\sum_{j\geq1}\floor{k/p^j}
       <\frac{k}{p-1}\leq k\leq2^r
       \leq v_p(q^{2^r}).
\]
For $p>k$ the factorial valuation is zero. Unique factorization gives the
divisibility. The $r$ individual squarings, not an exponentiation gate,
compute the output. The positive integer quotient by $k!$ need not be
computed.
\end{proof}

\begin{lemma}[Recording rational failures]\label{lem:compiler}
Let a rational straight-line program of length $s$ start with $1,2$ and
use binary $+,-,\times,/$, with nonzero rational divisors and nonzero
designated output. There is a nonzero integer $b$, computed using at most
$5s$ allowed integer gates, such that for every prime $p$, $p\mid b$ if
and only if the same rational program over $\F_p$ encounters division by
zero or completes with zero output. Instructions whose values are unused
later are retained as part of the program.
\end{lemma}
\begin{proof}
Represent every rational register $q_i$ by an unreduced pair of integer
registers $(U_i,V_i)$ with $q_i=U_i/V_i$ and $V_i\ne0$. Initially use
$(1,1),(2,1)$. The simulation rules are
\[
\begin{array}{c|cc}
q_i&U_i&V_i\\\hline
q_a+q_b&U_aV_b+U_bV_a&V_aV_b\\
q_a-q_b&U_aV_b-U_bV_a&V_aV_b\\
q_aq_b&U_aU_b&V_aV_b\\
q_a/q_b&U_aV_b&V_aU_b
\end{array}
\]
They cost at most four allowed gates per rational gate. All $V_i$ are
nonzero over $\Z$ by induction; at a division the numerator $U_b$ is
nonzero because $q_b\ne0$ over $\Q$. If $o$ is the output index, set
\[
 b=U_o\prod_{i:\,\text{division at }i}U_{b(i)},
\]
where $b(i)$ is that division's divisor index. Every factor is nonzero.
There are at most $s$ appended multiplications, giving the bound $5s$.

As long as the modular rational execution succeeds, induction gives
$V_i\not\equiv0\pmod p$ and modular value $U_i/V_i$. At its first failed
division, the corresponding $U_{b(i)}$ vanishes, so $p\mid b$. If execution
completes with zero output, $U_o$ vanishes. If it completes with nonzero
output, every recorded divisor numerator and $U_o$ is nonzero, so $p\nmid b$.
This also shows why unused test inverses must not be optimized away.
\end{proof}

\section{A family with a non-torsion marked class}\label{sec:family}
Let $H\geq1$ and $z=(z_1,\ldots,z_5)\in\{0,\ldots,H-1\}^5$. Set
\begin{equation}\label{eq:parameters}
 t_1=1+3z_1,\qquad t_i=3z_i\quad(2\leq i\leq5),
\end{equation}
and define
\begin{equation}\label{eq:family}
 f_z(X)=(1+3^{10}t_5)X^5+
             \sum_{i=1}^4 3^{2i}t_iX^i-\sum_{i=1}^5t_i.
\end{equation}
Write $a=1+3^{10}t_5$. In characteristic zero let $C_z$ be the smooth
completion of $Y^2=f_z(X)$, let $\infty$ be its point at infinity, and put
\begin{equation}\label{eq:point}
 P=(3^{-2},3^{-5}),\qquad D=[P-\infty]\in\Jac(C_z)(\Q).
\end{equation}

\begin{lemma}\label{lem:infinite}
For every permitted integer vector $z$, $C_z$ is a smooth genus-two curve
with good reduction at $3$, and $D$ has infinite order.
\end{lemma}
\begin{proof}
Substitution gives $f_z(3^{-2})=3^{-10}$, proving that $P$ lies on the curve.
Modulo $3$, $f_z$ is $X^5-1$, with derivative $2X^4$ and no common root.
The leading coefficient is a $3$-adic unit. Thus the quintic is squarefree
over $\Q$ and defines a smooth proper genus-two curve over $\Z_3$.

For later arithmetic pass to the monic integral model $y^2=G(x)$ by
$x=aX$, $y=a^2Y$. If $f_z(X)=\sum_{i=0}^5a_iX^i$, then
\begin{equation}\label{eq:monic}
 G(x)=x^5+a_4x^4+a_3ax^3+a_2a^2x^2+a_1a^3x+a_0a^4.
\end{equation}
The transformed point is $(a/9,a^2/243)$. At infinity, $t=x^2/y$ is
an integral local parameter and its value at this point is $3$. More
explicitly, with $v=1/x$ the infinity chart satisfies
$t^2=v/(1+a_4v+a_3av^2+\cdots+a_0a^4v^5)$; its denominator is a unit
at $v=0$. At the point, $v=9/a\in9\Z_3$. Thus $P$ reduces to infinity,
and $D$ lies in the kernel of reduction of the Jacobian.

The class $D$ is nonzero: otherwise a function with divisor $P-\infty$
would give a degree-one morphism to $\mathbb P^1$, impossible for a
genus-two curve. We prove that the kernel of reduction is torsion-free
in the present situation. Put $J=\Jac(C_z)$, let $\widetilde J$ be its
reduction modulo $3$, and write
\[
 J_1(\Q_3)=\ker\bigl(J(\Q_3)\longrightarrow\widetilde J(\F_3)\bigr).
\]
The monic model above has coefficients in $\Z_3$, unit discriminant,
and a rational point at infinity. Applying~\cite[Corollary~4.5 and
Theorem~4.6]{Grant} with $R=\Z_3$ and maximal ideal
$3\Z_3$, integral local parameters $\xi_1,\xi_2$ at the
identity $O$ give a bijection
\[
 \theta:J_1(\Q_3)\longrightarrow(3\Z_3)^2,\qquad
 Q\longmapsto\bigl(\xi_1(Q),\xi_2(Q)\bigr),\qquad \theta(O)=(0,0).
\]
Under this bijection the group operation is given by a commutative
formal group law $F$ over $\Z_3$. In particular,
$\theta(D)\in(3\Z_3)^2\setminus\{(0,0)\}$, since $D$ is nonzero and
reduces to $O$. More generally, every nonzero point in $J_1(\Q_3)$ has
a nonzero coordinate vector $s=(s_1,s_2)$ in this product, so
$\min_i v_3(s_i)\geq1$, with the convention $v_3(0)=+\infty$.
The formal group series converge on these coordinates.

Write $F(s,t)=(F_1(s,t),F_2(s,t))$, with
$s=(s_1,s_2)^{\top}$ and $t=(t_1,t_2)^{\top}$ regarded as column
vectors. For $m\geq1$, let $\Phi_m(s)=[m](s)$ denote the formal power
series for multiplication by $m$ in this group, and define the
$2$-by-$2$ matrices
\[
 J_{\Phi_m}(s)=
 \left(\frac{\partial(\Phi_m)_i}{\partial s_j}(s)\right)_{i,j=1}^2,
 \qquad
 A(s)=
 \left(\frac{\partial F_i}{\partial t_j}(s,0)\right)_{i,j=1}^2.
\]
Here $J_{\Phi_m}$ is the matrix of formal partial derivatives, unrelated
to the divisor class $D$. Since $F(s,t)=s+t+$ terms of total degree at
least two, $A(0)=I_2$, while $\Phi_m(s)=ms+$ terms of degree at least
two gives $J_{\Phi_m}(0)=mI_2$. In particular, $A(s)$ is invertible over
$\Z_3[[s_1,s_2]]$.
Commutativity makes multiplication by $m$ a group homomorphism, so
\[
 \Phi_m(F(s,t))=F(\Phi_m(s),\Phi_m(t)).
\]
Differentiating this identity with respect to $t$ and then setting $t=0$,
the formal chain rule, together with $F(s,0)=s$ and $\Phi_m(0)=0$, yields
\[
 J_{\Phi_m}(s)A(s)
 =A(\Phi_m(s))J_{\Phi_m}(0)
 =mA(\Phi_m(s)).
\]
For $m=3$, multiplication on the right by $A(s)^{-1}$ shows that
$J_{\Phi_3}(s)\equiv0\pmod3$. Thus every formal partial derivative of
each component of $\Phi_3$ vanishes modulo $3$, and every monomial with
nonzero coefficient in $\Phi_3\bmod3$ has each variable exponent
divisible by $3$.
In particular the degree-two coefficients are divisible by $3$.
The linear term is $3s$. For a nonzero point of $J_1(\Q_3)$, let $s$
be its coordinates under $\theta$ and put $v=\min_i v_3(s_i)$.
As explained above, $1\leq v<+\infty$. The linear term has minimum
valuation $v+1$;
degree-two terms have valuation at least $2v+1\geq v+2$, and terms of
degree at least three have valuation at least $3v\geq v+2$.
Thus $[3](s)\ne0$ and its minimum valuation is $v+1$.
For $3\nmid m$, the linear term $ms$ dominates the higher terms, so
$[m](s)$ has minimum valuation $v$. Factoring any positive integer as
$3^rm$ proves torsion-freeness of $J_1(\Q_3)$. Since $D$ is a nonzero
point of this kernel, it has infinite order.
\end{proof}

\section{A charged nonzero integer certificate}\label{sec:arithmetic}
\begin{proposition}\label{prop:certificate}
There is an absolute constant $C_0$ such that for every $N,H\geq1$ and
every permitted vector $z$ there exists a positive integer $Z(N,z)$ with
\begin{equation}\label{eq:cert-cost}
 \tau(Z(N,z))\leq C_0\bigl(1+\ceil{\log_2N}+\ceil{\log_2(H+1)}\bigr).
\end{equation}
It is divisible by $2,3$. For every prime $p\geq5$, it is divisible by
$p$ whenever $f_z$ modulo $p$ fails to be a squarefree quintic, or when
$f_z$ is a squarefree quintic and $[N]D_p=0$ in $\Jac(C_z)(\F_p)$.
The circuit is fixed from the rational data $N,z$ independently of $p$.
\end{proposition}
\begin{proof}
We first specify complete group arithmetic rather than using a generic
addition formula. Over any field of characteristic different from two,
each class on the monic curve~\eqref{eq:monic} has a unique reduced
Mumford pair $(u,v)$, with $u$ monic,
\[
 \deg v<\deg u\leq2,\qquad u\mid G-v^2;
\]
the identity is $(1,0)$. Given two reduced pairs, compute a monic gcd and
Bezout coefficients
\[
 d=\gcd(u_1,u_2,v_1+v_2)
   =e_1u_1+e_2u_2+e_3(v_1+v_2).
\]
Cantor composition~\cite[Theorem~10.3.14]{Galbraith} gives
\begin{equation}\label{eq:cantor}
 u=\frac{u_1u_2}{d^2},\qquad
 v=\frac{e_1u_1v_2+e_2u_2v_1+e_3(v_1v_2+G)}{d}\pmod u.
\end{equation}
These polynomial divisions are exact. If $\deg u>2$, replace the pair by
\begin{equation}\label{eq:reduction}
 u'=\operatorname{monic}((G-v^2)/u),\qquad v'=-v\pmod{u'}.
\end{equation}
The cited composition theorem includes coincident and opposite support;
no coprimality of $u_1,u_2$ is assumed. The reduction is
\cite[\S2, Algorithm~1]{KS}, and uniqueness is
\cite[Theorem~10.4.1]{Galbraith}. The composed $u$ has degree at most
four and the remainder $v$ has degree at most three. If $\deg u=4$,
the next degree is at most two. If $\deg u=3$, then $\deg v\leq2$ and
the next degree is two. Hence one reduction suffices. Degree-zero identity
cases can be handled directly.

Use coefficient arrays, schoolbook multiplication, long division and two
extended Euclidean algorithms to implement these formulas. Every degree
is bounded by an absolute constant. Each nonzero Euclidean remainder
strictly decreases degree. There are consequently an absolute bounded
number of scalar $+,-,\times,/$ operations and coefficient zero tests per
addition, including all exceptional branches. Zero inputs to a gcd are
handled by their usual finite branches. Thus there is a universal bound
$M$ for both scalar gates and tested coefficients in one addition.

Compute $[N]D$ over $\Q$ by the binary double-and-add schedule, which
uses at most $2\floor{\log_2N}$ additions. For the particular literal
$N,z$, fix its exact rational branch path. At every coefficient test
whose rational value $c$ is nonzero, append and retain the scalar
instruction $1/c$, even if its value is never used. Rationally zero tests
require no inverse. This includes every degree, leading-coefficient and
normalization test; the inverse tests can be inserted before the relevant
branch arithmetic. Lemma~\ref{lem:infinite} ensures that every nonzero
scalar prefix, and in particular $[N]D$, is nonzero over $\Q$.

At the start of the rational program construct all parameters and monic
coefficients, and retain $1/2,1/3,1/a,1/\Disc(G)$. Each literal $z_i<H$
can be constructed by its binary expansion using $O(1+\log(H+1))$
integer gates. The powers of $3$ and the polynomial coefficient formulas
have constant cost beyond those literal constructions. The discriminant
can be computed by expanding the fixed $9$-by-$9$ Sylvester determinant
of $G,G'$ into signed permutation products; this is a constant-size
coefficient circuit, not a determinant gate. Its nonvanishing over $\Q$
follows from good reduction at $3$. Initialize the reduced marked pair
$u=x-a/9$, $v=a^2/243$. Append the designated nonzero output $1$.
All divisors in the fixed rational program are nonzero over $\Q$ and its
length is
\begin{equation}\label{eq:rational-length}
 s\leq C_1(1+\ceil{\log_2N}+\ceil{\log_2(H+1)})
\end{equation}
for an absolute constant $C_1$. The inserted inverses are included in it.

Consider a prime $p\geq5$. If $a=0$ modulo $p$, the retained $1/a$
fails, unless some earlier division has already failed. If $a$ is a unit
but the quintic is not squarefree, the coordinate change is invertible
and $\Disc(G)=0$, so its retained inverse fails. Now suppose the quintic
is squarefree and the entire rational program succeeds modulo $p$.
Inductively, every rationally zero test specializes to zero, and each
rationally nonzero test specializes to a nonzero value because its
retained inverse succeeds. All degrees, Euclidean steps and branches
are therefore unchanged. The formulas~\eqref{eq:cantor}--\eqref{eq:reduction}
then specialize to correct reduced representatives over $\F_p$.
The final $u$ has positive degree over $\Q$, by infinite order of $D$,
and the same degree modulo $p$. It cannot represent the identity $(1,0)$.
Consequently $[N]D_p=0$ forces a modular division failure.

Apply Lemma~\ref{lem:compiler} to this entire rational program and square
its nonzero integer output $b$. Then $Z(N,z)=b^2>0$ has all the required
prime divisors and costs at most $5s+1$ gates. The retained inverses of
$2,3$ supply those prime factors. This proves~\eqref{eq:cert-cost}.
No branch, zero-test, gcd, division or group-operation gate remains in
the integer program. All rational branches were selected before the
circuit was fixed, without reference to a tested prime.
\end{proof}

\section{Smooth orders on the actual common source}\label{sec:density}
\begin{lemma}[Transfer to a marked-point slice]\label{lem:slice}
Fix a prime $p>8100$ and a set $\mathcal A\subseteq I_p\cap\Z$.
If $z$ is uniform in $\F_p^5$ and $f_z$ is defined by
\eqref{eq:parameters}--\eqref{eq:family}, then
\begin{equation}\label{eq:slice-rate}
 \Pr\bigl(f_z\text{ is a squarefree quintic and }
          \#\Jac(C_z)(\F_p)\in\mathcal A\bigr)
 \geq \frac{\max(|\mathcal A|-p,0)}
 {48000p^{3/2}(\log p)^2(\log\log p)^2}.
\end{equation}
\end{lemma}
\begin{proof}
Let $\mathcal S$ be all squarefree quintics with Jacobian order in
$\mathcal A$. Theorem~\ref{thm:lpp2} gives
\[
 |\mathcal S|\geq
 \frac{\max(|\mathcal A|-p,0)p^{9/2}}
      {24000(\log p)^2(\log\log p)^2}.
\]
Let $\mathcal T$ consist of those squarefree quintics $g$ with $g(0)=1$
and order in $\mathcal A$. For $f\in\mathcal S$ and an affine point
$(c,d)$ with $d\ne0$ define $g(W)=f(W+c)/d^2$. Conversely, for
$g\in\mathcal T$, $c\in\F_p$, $d\in\F_p^*$, set
$f(X)=d^2g(X-c)$. This is a bijection of the corresponding triples,
preserving curve isomorphism and Jacobian order. Thus
\[
 |\mathcal T|p(p-1)
 =\sum_{f\in\mathcal S}\#\{(c,d):d^2=f(c),\ d\ne0\}.
\]
By~\eqref{eq:weil}, each summand is at least $p-4\sqrt p-5$:
remove the point at infinity and at most five zero-ordinate points.
Therefore
\[
 \frac{|\mathcal T|}{p^5}\geq
 \frac{|\mathcal S|}{p^6}
 \frac{p-4\sqrt p-5}{p-1}
 \geq\frac{|\mathcal S|}{2p^6}.
\]
The last inequality holds for $p>8100$.

Put $c_0=3^{-2}$, $d_0=3^{-5}$ in $\F_p$. The five affine maps from
$z$ to the coefficients of $X,\ldots,X^5$ in $f_z$ are bijections,
and the constant coefficient is exactly the one imposing
$f_z(c_0)=d_0^2$. Hence uniform $z$ samples the entire five-dimensional
slice of polynomials of degree at most five with that prescribed point.
The map $g(W)=f(W+c_0)/d_0^2$ identifies this slice with all polynomials
with constant term one, a set of size $p^5$. It preserves the required
properties. Thus the desired probability equals $|\mathcal T|/p^5$,
which proves the claim. Leading coefficients have not been fixed.
\end{proof}

For $k\geq2$ define the positive real quantities
\begin{equation}\label{eq:ASE}
 A=\log(16k^2),\qquad
 S=A^{2/3}(\log A)^{1/3},\qquad
 E=A^{1/3}(\log A)^{8/3},
\end{equation}
and put
\begin{equation}\label{eq:BdeltaH}
 B=\ceil{\exp(1000S)},\qquad
 \delta=\frac{\exp(-E)}{192000A^4},\qquad H=10k.
\end{equation}

\begin{proposition}\label{prop:density}
There are absolute integers $p_0>8100$ and $K_0\geq p_0$ such that for
every $k\geq K_0$, the single sampler with independent coordinates
$z_i$ uniform in $\{0,\ldots,H-1\}$ satisfies, for every prime
$p_0\leq p\leq k$,
\begin{equation}\label{eq:actual-density}
 \Pr\bigl(f_z\text{ is a squarefree quintic mod }p\text{ and }
          \#\Jac(C_z)(\F_p)\text{ is }B\text{-smooth}\bigr)\geq\delta.
\end{equation}
\end{proposition}
\begin{proof}
Set $x_p=p^2-\tfrac12p^{3/2}$, $y_p=p^{3/2}$ and
$e_p=(\log x_p)^{1/3}(\log\log x_p)^{8/3}$.
For all sufficiently large $p$, $x_p\geq x_*$ and
$x_p^{3/5}\leq y_p\leq x_p$. Theorem~\ref{thm:lpp1} supplies at least
$p^{3/2}e^{-e_p}$ smooth integers in $I_p$. Also $e_p=o(\log p)$, so
$p^{3/2}e^{-e_p}\geq2p$ for every sufficiently large $p$.
Choose a fixed $p_0>8100$ so these statements all hold for every $p\geq p_0$.

If $p\leq k$, then $\log x_p\leq A$. Both
$t^{2/3}(\log t)^{1/3}$ and $t^{1/3}(\log t)^{8/3}$ are increasing for
$t>1$, and $c_*<1000$. Hence their smoothness threshold is at most $B$
and $e_p\leq E$. If $\mathcal A$ is the set of $B$-smooth integers in
$I_p$, then
\[
 |\mathcal A|-p\geq\tfrac12p^{3/2}e^{-e_p}
                        \geq\tfrac12p^{3/2}e^{-E}.
\]
Deleting every exceptional order is therefore harmless at this stated loss.
Lemma~\ref{lem:slice} gives uniform-residue probability at least
\[
 \frac{e^{-E}}{96000(\log p)^2(\log\log p)^2}
 \geq\frac{e^{-E}}{96000A^4}.
\]

For the actual integer sampler, each residue in one coordinate has
probability at least $\floor{H/p}/H\geq(1-p/H)/p\geq(9/10)/p$.
Coordinate independence gives every residue vector at least $(9/10)^5$
times its uniform mass. Since $(9/10)^5>1/2$, its event probability is
at least $\delta$. All choices are independent of the tested prime.
\end{proof}

\section{A common scalar and simultaneous prime coverage}\label{sec:cover}
\begin{lemma}\label{lem:group-size}
For a squarefree quintic over $\F_p$, $p\geq5$,
$\#\Jac(C)(\F_p)\leq16p^2$.
\end{lemma}
\begin{proof}
Pass to the monic model and count reduced pairs. Degree zero contributes
one. Degree one contributes at most $2p$. There are $p^2$ monic quadratic
$u$, each admitting at most four possible $v$ with $v^2=G\bmod u$.
For irreducible $u$ the quotient is a field, giving at most two roots;
for two distinct linear factors the Chinese remainder theorem gives at
most four. For $u=(x-c)^2$, write $v=v_0+v_1(x-c)$. The equations are
$v_0^2=G(c)$ and $2v_0v_1=G'(c)$. If $G(c)\ne0$ they give at most
two solutions; if $G(c)=0$, squarefreeness gives $G'(c)\ne0$, so there
are none. Thus the total is at most $1+2p+4p^2\leq16p^2$.
\end{proof}

Using~\eqref{eq:ASE}--\eqref{eq:BdeltaH}, choose the positive integer
\begin{equation}\label{eq:N}
 N=\prod_{\substack{\ell\leq B\\\ell\text{ prime}}}
           \ell^{\floor{A/\log\ell}}.
\end{equation}
Every $B$-smooth integer $n\leq16k^2$ divides $N$. Moreover,
\begin{equation}\label{eq:logN}
 \log N\leq A\#\{\ell\leq B:\ell\text{ prime}\}\leq AB.
\end{equation}
This definition specifies the binary schedule for scalar multiplication.
It supplies no free prime list, scalar register, logarithm or floor gate.

For each permitted literal $z$, fix one of the circuits in
Proposition~\ref{prop:certificate} from its rational data, independently of
any tested prime. On the event in~\eqref{eq:actual-density}, the Jacobian
order is at most $16p^2\leq16k^2$, so it divides $N$. Lagrange's theorem
gives $[N]D_p=0$. Therefore, for every $p_0\leq p\leq k$,
\begin{equation}\label{eq:capture}
 \Pr_z(p\mid Z(N,z))\geq\delta.
\end{equation}

Take
\begin{equation}\label{eq:J}
 J=\ceil{(\log k+1)/\delta}
\end{equation}
independent samples from the same integer grid. For each fixed tested prime,
the probability of missing it in all $J$ outputs is at most
$(1-\delta)^J\leq e^{-\delta J}\leq e^{-1}/k$.
There are fewer than $k$ such primes, so a union bound is strictly less
than one. Thus there exist literal vectors $z^{(1)},\ldots,z^{(J)}$
whose outputs jointly cover every prime in $[p_0,k]$. Independence
between different primes is neither asserted nor needed.

Construct $p_0!$ by successive integer increments and accumulator
multiplications, at cost at most $2p_0$. The positive integer
\begin{equation}\label{eq:Q}
 Q=p_0!\prod_{j=1}^J Z(N,z^{(j)})
\end{equation}
is divisible by every prime at most $k$. Concatenating the $J$ circuits
and multiplying their outputs into the accumulator gives
\begin{equation}\label{eq:total}
 \tau(Q)\leq2p_0+
 J\left[C_0\bigl(1+\ceil{\log_2N}+\ceil{\log_2(10k+1)}\bigr)+1\right].
\end{equation}
All parameters, tests and output multiplications are included.

\section{Asymptotic accounting and conclusion}\label{sec:asymptotics}
\begin{proof}[Proof of Theorem~\ref{thm:main}]
As $k\to\infty$, $A=2\log k+\log16$ and
\begin{equation}\label{eq:EoverS}
 \frac{E}{S}=A^{-1/3}(\log A)^{7/3}\longrightarrow0,
 \qquad \log A=o(S).
\end{equation}
Since $B\leq2e^{1000S}$, equation~\eqref{eq:logN} implies
$\ceil{\log_2N}\leq1+AB/\log2$. Also
$\ceil{\log_2(10k+1)}=O(A)$. Consequently the bracket in~\eqref{eq:total}
is $\exp(O(S))$. On the other hand,
\[
 J\leq1+192000(\log k+1)A^4e^E,
 \qquad \log J\leq E+O(\log A)=o(S).
\]
It follows that $\tau(Q)\leq\exp(O(S))$.

Apply Lemma~\ref{lem:amplify}. Appending $\ceil{\log_2k}$ ordinary
squarings produces $Q^{2^{\ceil{\log_2k}}}=mk!$ with $m>0$, at total
cost $\exp(O(S))$. This proves~\eqref{eq:allk} for all sufficiently
large $k$. For each of the finitely many remaining $k\geq2$, construct
$k!$ directly in at most $2k$ gates. The exponent scale $S$ defined
in~\eqref{eq:ASE} is positive throughout this range, so one absolute
increase of $C$ covers all those cases.

Finally $S\sim2^{2/3}(\log k)^{2/3}(\log\log k)^{1/3}$ and
$S/\log k\to0$, giving~\eqref{eq:main} and the stated consequence for
every fixed $\varepsilon>0$.
\end{proof}

\begin{remark}[Scope of the bound]
The bounding function $\exp(C(\log k)^{2/3}(\log\log k)^{1/3})$ exceeds
every fixed power of $\log k$ for large $k$; this observation concerns the
bound, not a lower bound for $\mu(k)$. The argument neither proves nor
refutes Smale's original second or third question. It gives no uniform
bit-complexity bound for finding the circuit. In particular, bounding the
number of arithmetic gates does not bound the bit lengths of intermediate
integer registers.
\end{remark}

\section*{Acknowledgements}
The author thanks Xiaohua (Michael) Xuan for introducing him to this research
topic and for his continued support.


\begin{thebibliography}{99}
\bibitem{AGHT}
Alekseev, Y., Grigoriev, D., Hirsch, E. and Tzameret, I.,
\emph{\href{https://doi.org/10.1137/20M1374523}{Semialgebraic proofs,
IPS lower bounds, and the $\tau$-conjecture: Can a natural number be negative?}},
SIAM J. Comput. \textbf{53} (2024), no.~3, 648--700.

\bibitem{Burgisser}
B\"urgisser, P., \emph{\href{https://doi.org/10.1007/s00037-009-0260-x}
{On defining integers and proving arithmetic circuit lower bounds}},
Comput. Complexity \textbf{18} (2009), no.~1, 81--103.

\bibitem{Cantor}
Cantor, D.~G., \emph{Computing in the Jacobian of a hyperelliptic curve},
Math. Comp. \textbf{48} (1987), no.~177, 95--101.

\bibitem{Cheng}
Cheng, Q., \emph{\href{https://doi.org/10.1016/j.tcs.2004.06.020}
{On the ultimate complexity of factorials}},
Theoret. Comput. Sci. \textbf{326} (2004), nos.~1--3, 419--429.
\href{https://qcheng2023.github.io/paper/factorial.pdf}
{Author preprint, May~25, 2004}.

\bibitem{dMS}
de Melo, W. and Svaiter, B.~F.,
\emph{\href{https://doi.org/10.1090/S0002-9939-96-03173-5}
{The cost of computing integers}},
Proc. Amer. Math. Soc. \textbf{124} (1996), no.~5, 1377--1378.

\bibitem{Galbraith}
Galbraith, S.~D., \emph{Mathematics of Public Key Cryptography},
Cambridge University Press, Cambridge, 2012;
\href{https://www.math.auckland.ac.nz/~sgal018/crypto-book/ch10.pdf}
{online extended and corrected version~2.0, Chapter~10}.

\bibitem{Grant}
Grant, D., \emph{\href{https://doi.org/10.1515/crll.1990.411.96}
{Formal groups in genus two}},
J. reine angew. Math. \textbf{411} (1990), 96--121.

\bibitem{KS}
Kohel, D.~R. and Smith, B.~A.,
\emph{Efficiently computable endomorphisms for hyperelliptic curves},
\href{https://arxiv.org/abs/math/0603505v1}{arXiv: math/0603505v1}, 2006.

\bibitem{Koiran}
Koiran, P., \emph{\href{https://doi.org/10.1007/s00037-004-0186-2}
{Valiant's model and the cost of computing integers}},
Comput. Complexity \textbf{13} (2004), nos.~3--4, 131--146.

\bibitem{LPP1}
Lenstra, H.~W., Jr., Pila, J. and Pomerance, C.,
\emph{\href{https://doi.org/10.1098/rsta.1993.0138}
{A hyperelliptic smoothness test. I}},
Philos. Trans. Roy. Soc. London Ser. A \textbf{345} (1993), 397--408.

\bibitem{LPP2}
Lenstra, H.~W., Jr., Pila, J. and Pomerance, C.,
\emph{\href{https://doi.org/10.1112/plms/84.1.105}
{A hyperelliptic smoothness test, II}},
Proc. London Math. Soc. (3) \textbf{84} (2002), 105--146.

\bibitem{Markstrom}
Markstr\"om, K., \emph{\href{https://math.colgate.edu/~integers/p6/p6.pdf}
{The straight line complexity of small factorials and primorials}},
Integers \textbf{15} (2015), Paper No.~A6, 8~pp.

\bibitem{Milne}
Milne, J.~S., \emph{Jacobian varieties},
in \emph{Arithmetic Geometry} (Storrs, Conn., 1984),
167--212, Springer, New York, 1986;
\href{https://www.jmilne.org/math/xnotes/JVs.pdf}
{revised author notes, June~12, 2021}.

\bibitem{Moreira}
Moreira, C.~G.~T.~de~A.,
\emph{\href{https://doi.org/10.1090/S0002-9939-97-03583-1}
{On asymptotic estimates for arithmetic cost functions}},
Proc. Amer. Math. Soc. \textbf{125} (1997), no.~2, 347--353.

\bibitem{SS}
Shub, M. and Smale, S., \emph{On the intractability of Hilbert's
Nullstellensatz and an algebraic version of ``$NP\ne P$?''},
Duke Math. J. \textbf{81} (1995), no.~1, 47--54.

\bibitem{Smale}
Smale, S., \emph{\href{https://doi.org/10.1007/BF03025291}
{Mathematical problems for the next century}},
Math. Intelligencer \textbf{20} (1998), no.~2, 7--15.
\end{thebibliography}
\end{document}